\pdfoutput=1
\RequirePackage{ifpdf}
\ifpdf % We are running pdfTeX in pdf mode
\documentclass[pdftex]{sigma}
\else
\documentclass{sigma}
\fi
\begin{document}

\allowdisplaybreaks

\renewcommand{\thefootnote}{$\star$}

\renewcommand{\PaperNumber}{027}

\FirstPageHeading

\ShortArticleName{$G$-Strands and Peakon Collisions on $\operatorname{Dif\/f}(\mathbb{R})$}

\ArticleName{$\boldsymbol{G}$-Strands and Peakon Collisions on $\boldsymbol{\operatorname{Dif\/f}(\mathbb{R})}\,$\footnote{This paper
is a~contribution to the Special Issue ``Symmetries of Dif\/ferential Equations: Frames, Invariants and~Applications''.
The full collection is available at
\href{http://www.emis.de/journals/SIGMA/SDE2012.html}{http://www.emis.de/journals/SIGMA/SDE2012.html}}}

\Author{Darryl D.~HOLM~$^\dag$ and Rossen I.~IVANOV~$^\ddag$}

\AuthorNameForHeading{D.D.~Holm and R.I.~Ivanov}

\Address{$^\dag$~Department of Mathematics, Imperial College London, London SW7 2AZ, UK}
\EmailD{\href{mailto:d.holm@ic.ac.uk}{d.holm@ic.ac.uk}}
\URLaddressD{\url{http://www2.imperial.ac.uk/~dholm/}}

\Address{$^\ddag$~School of Mathematical Sciences, Dublin Institute of Technology,
\\
\hphantom{$^\ddag$}~Kevin Street, Dublin 8, Ireland}
\EmailD{\href{mailto:rivanov@dit.ie}{rivanov@dit.ie}}

\ArticleDates{Received October 29, 2012, in f\/inal form March 21, 2013; Published online March 26, 2013}

\Abstract{A $G$-strand is a~map $g:\mathbb{R}\times\mathbb{R}\to G$ for a~Lie group $G$ that follows from Hamilton's
principle for a~certain class of $G$-invariant Lagrangians.
Some $G$-strands on f\/inite-dimensional groups satisfy $1+1$ space-time evolutionary equations that admit soliton
solutions as completely integrable Hamiltonian systems.
For example, the ${\rm SO}(3)$-strand equations may be regarded physically as integrable dynamics for solitons on
a~continuous spin chain.
Previous work has shown that $G$-strands for dif\/feomorphisms on the real line possess solutions with singular support
(e.g.~peakons).
This paper studies \emph{collisions} of such singular solutions of $G$-strands when $G={\rm
Dif\/f}(\mathbb{R})$ is the group of dif\/feomorphisms of the real line $\mathbb{R}$, for which the group product is
composition of smooth invertible functions.
In the case of peakon-antipeakon collisions, the solution reduces to solving either Laplace's equation or the wave
equation (depending on a~sign in the Lagrangian) and is written in terms of their solutions.
We also consider the complexif\/ied systems of $G$-strand equations for $G={\rm Dif\/f}(\mathbb{R})$ corresponding to
a~harmonic map $g: \mathbb{C}\to{\rm Dif\/f}(\mathbb{R})$ and f\/ind explicit expressions for its peakon-antipeakon
solutions, as well.
}


\Classification{37J15; 37K05; 35R01}
\begin{flushright}
\begin{minipage}
{72mm}\it In honor of Peter Olver's sixtieth birthday.
\\
Happy birthday, Peter!
\end{minipage}
\end{flushright}

\renewcommand{\thefootnote}{\arabic{footnote}} \setcounter{footnote}{0}

\section[Euler-Poincar\'e equations for a~$G$-strand]{Euler--Poincar\'e equations for a~$\boldsymbol{G}$-strand}
\label{intro-sec}

The Euler--Poincar\'e (EP) theory of $G$-strands is an extension of the classical chiral models to include norms that
are not bi-invariant~\cite{Ho2011,HoIvPe2012}.
The classical chiral models are nonlinear relativistically invariant Lagrangian f\/ield theories on group manifolds.
As such, they are fundamental in theoretical physics.
The vast literature of results for these models is fascinating.
For example, it is well known that these models are integrable in $1 + 1$ dimensions and possess soliton solutions.
See~\cite{Ch1981,deVega1979,Mai86,Man1994,NoMaPiZa1984,Po1976,ZaMi1980,ZaMi1978} and references therein for
discussions of the many aspects of integrability of the chiral models, including the famous dressing method for
explicitly deriving the soliton solutions of these models, which is given in~\cite{ZaMi1980}.
The solitons for the ${\rm O}(3)$ chiral model are particularly familiar~\cite{Ya1988}, because this model allows an
integrable reduction to the one-component sine-Gordon equation (e.g.\
see~\cite{NoMaPiZa1984}).
Many generalizations of these models have been introduced.
For example, generalized chiral models with metrics that are not ad-invariant on the Lie algebra are considered
in~\cite{So1997}.
Other generalizations of chiral models for non-semisimple groups are studied in~\cite{HlaSno2001}.
An integrable chiral model in $2+1$ dimensions was proposed in~\cite{Wa1988} and the possibility for fermionic
interpretation of the current variables was explored in~\cite{Wi1984}.
Finally, their formulation as Euler--Poincar\'e equations for $G$-strands with Lag\-rangians given as norms on tangent
spaces of f\/inite-dimensional Lie groups was accomplished in~\cite{Ho2011,HoIvPe2012}.

The class of $G$-strand problems treated here generalizes the classical chiral model, which appear as a~special case of
the $G$-strands when their structure group is ${\rm SO}(3)$ and the norm def\/ining their Lagrangian is bi-invariant, as
discussed in~\cite{HoIvPe2012}.
$G$-strands are found here to also be associated with a~certain generalization of harmonic maps into Lie groups.
Harmonic maps into Lie groups were studied mathematically, and their relation to the chiral model via
complexif\/ication of its independent coordinates was discussed in~\cite{Uh1989}.
Following this same idea, we shall complexi\-fy the \emph{independent coordinates} of the $G$-strand on the Lie group of
dif\/feomorphisms ${\rm Dif\/f}(\mathbb{R})$ and use the Euler--Poincar\'e theory to derive a~class of equations which
may be regarded as partial dif\/ferential equations for harmonic maps of $\mathbb{C}$ into ${\rm Dif\/f}(\mathbb{R})$.
This particular advance is not by itself the main point of the present paper, though.
The main point is the investigation of the interactions between singular solutions of the $G$-strand equations for the
group of dif\/feomorphisms ${\rm Dif\/f}(\mathbb{R})$, both as evolutionary equations and as harmonic maps.

{\bf Main content of the paper.}
 Section~\ref{intro-sec} summarizes the derivation of the class of $G$-strand equations f\/irst studied in the
context of Euler--Poincar\'e theory in~\cite{Ho2011,HoIvPe2012}.
  Section~\ref{DiffStrand-sec} discusses singular (peakon) solutions and their collision dynamics for the
$G$-strand equations with $G={\rm Dif\/f}(\mathbb{R})$.
These peakon solutions were shown to exist in~\cite{HoIvPe2012}.
Here we study the dynamics of their pairwise interactions.
  Section~\ref{complex-sec} introduces the complexif\/ied ${\rm Dif\/f}(\mathbb{R})$-strand equations and
determines their solutions corresponding to peakon collision dynamics.
 Section~\ref{conclusion-sec} summarizes the results of the paper and provides some outlook for future research.

{\bf Left $\boldsymbol{G}$-invariant Lagrangian.} We begin with the following ingredients of EP theory.
For more details and discussion, see~\cite{Ho2011,HoMaRa1998}.
\begin{itemize}\itemsep=0pt
\item Let $G$ be a~Lie group.
The map $g: (t,s)\in\mathbb{R}\times\mathbb{R}\to g(t,s)\in G$ possesses two types of tangent vectors, $(u_g,v_g)\in
TG\times TG$, corresponding to its two types of derivatives, with respect to $t$ and $s$.
\item Assume that the function $L(g,u_g,v_g):T G\times TG\rightarrow\mathbb{R}$ is right $G$-invariant.

\item Right $G$--invariance of $L$ permits us to def\/ine $l:\mathfrak{g}\times\mathfrak{g}\rightarrow\mathbb{R}$ by
\begin{gather*}
l\big(u_gg^{-1},v_g g^{-1}\big)=L(g,u_g,v_g).
\end{gather*}
Conversely, this relation def\/ines for any $l:\mathfrak{g}\times\mathfrak{g}\rightarrow\mathbb{R}$ a~right
$G$-invariant function $L:T G\times TG\rightarrow\mathbb{R}$.
\item For a~map $g(t,s):\mathbb{R}\times\mathbb{R}\to G$ one def\/ines the right $G$-invariant tangent vectors at the
identity of $G$, given by
\begin{gather*}
{u}(t,s):=g_t g^{-1}(t,s)
\qquad
\hbox{and}
\qquad
{v}(t,s):=g_sg^{-1}(t,s).
\end{gather*}
\end{itemize}

\begin{lemma}
The right-trivialized tangent vectors $u(t,s)$ and $v(t,s)$ at the identity of~$G$ satisfy
\begin{gather}
v_t-u_s={\rm ad}_{u}{v}.
\label{zero-curv}
\end{gather}
\end{lemma}

\begin{proof}
The proof is standard and follows from equality of cross derivatives $g_{ts}=g_{st}$, cf.~\cite{Ho2011,HoIvPe2012, HoMaRa1998}.
As a~consequence, equation~\eqref{zero-curv} is often called a~\emph{zero-curvature relation}.
\end{proof}
\begin{theorem}[Euler--Poincar\'e theorem]\label{lall}
With the preceding notation, the following four statements are equivalent~{\rm \cite{Ho2011,HoIvPe2012}}:
\begin{enumerate}\itemsep=0pt
\item [$(i)$] Hamilton's variational principle on $T G\times TG$
\begin{gather*}
%\label{hamiltonprinciple}
\delta\int_{t_1}^{t_2}L(g(t,s),{g_t}(t,s),g_s(t,s))\,dsdt=0
\end{gather*}
holds, for variations $\delta g(t,s)$ of $g(t,s)$ vanishing at the endpoints in $t$ and $s$.
\item [$(ii)$] The function $g(t,s)$ satisfies Euler--Lagrange equations for $L$ on $G$, given by
\begin{gather*}
%\label{EL-eqns}
\frac{\partial L}{\partial g}-\frac{\partial}{\partial t}\frac{\partial L}{\partial g_t}-\frac{\partial}{\partial s}\frac{\partial L}{\partial g_s}=0
\end{gather*}

\item [$(iii)$] The constrained variational principle\footnote{As with the basic Euler--Poincar\'e
equations~\cite{MaRa1999}, this is not strictly a~variational principle in the same sense as the standard Hamilton's
principle.
It is more like the Lagrange d'Alembert principle, because we impose the stated constraints on the variations
allowed~\cite{HoMaRa1998}.
The proof of equations~\eqref{epvariations} follows the same pattern as the proof of the f\/irst zero-curvature
relation~\eqref{zero-curv}.}
\begin{gather*}
%\label{variationalprinciple}
\delta\int_{t_1}^{t_2}l({u}(t,s),{v}(t,s))\,ds dt=0
\end{gather*}
holds on $\mathfrak{g}\times\mathfrak{g}$, using variations of ${u}$ and ${v}$ of the forms
\begin{gather}
\label{epvariations}
\delta{u}={w_t}-{\rm ad}_{u}{w}
\qquad
\hbox{and}
\qquad
\delta{v}={w_s}-{\rm ad}_{v}{w},
\end{gather}
where $w(t,s)\in\mathfrak{g}$ vanishes at the endpoints.
\item [$(iv)$] The Euler--Poincar\'{e} equations hold on $\mathfrak{g}^*\times\mathfrak{g}^*$
\begin{gather}
\frac{d}{dt}\frac{\delta l}{\delta{u}}+\operatorname{ad}_{{u}}^{\ast}\frac{\delta l}{\delta{u}}+\frac{d}{ds}\frac{\delta l}{\delta{v}}+\operatorname{ad}_{{v}}^{\ast}\frac{\delta l}{\delta{v}}=0.
\label{EPst-eqn}
\end{gather}
\end{enumerate}
\end{theorem}

{\bf Evolutionary $\boldsymbol{G}$-strand.} We now def\/ine the fundamental quantity of interest in the remainder of the paper.
\begin{definition}
A {\it $G$-strand} is an evolutionary map into a~Lie group $G$, $g(t,{s}):\,\mathbb{R}\times\mathbb{R}\to G$, whose
dynamics in $(t,{s})\in\mathbb{R}\times\mathbb{R}$ may be obtained from Hamilton's principle for a~$G$-invariant
reduced Lagrangian $l:\mathfrak{g}\times\mathfrak{g}\to\mathbb{R}$, where $\mathfrak{g}$ is the Lie algebra of the
group $G$.
The $G$-strand system of evolutionary partial dif\/ferential equations for a~right $G$-invariant reduced Lagrangian
consists of the \emph{zero-curvature} equation~\eqref{zero-curv} and the Euler--Poincar\'e (EP) variational
equations~\eqref{EPst-eqn}.\footnote{The left $G$-invariant case yields similar equations, but with the opposite signs of ${\rm ad}$ and ${\rm
ad}^*$, as in~\cite{HoIvPe2012}.}

\end{definition}

Subclasses of the $G$-strand maps contain the principal chiral models of f\/ield theory in theore\-ti\-cal physics,
reviewed, e.g., in~\cite{Man1994}.
An interpretation of the $G$-strand equations as the dynamics of a~continuous spin chain is given in~\cite{Ho2011}.
This is the source of the term, `strand'.
The corresponding theory of molecular strands is discussed in~\cite{ElGBHoPuRa2010}.
Recently, a~covariant f\/ield theory of $G$-strands in higher dimensions for $G={\rm Dif\/f}(\mathbb{R})$ has also been
developed~\cite{FGB2012}.

{\bf Overview and organization of the paper.} In Section~\ref{DiffStrand-sec}, we shall f\/irst summarize
$G$-strand dynamics for dif\/feomorphisms on the real line and recall from~\cite{HoIvPe2012} that the
Dif\/f$(\mathbb{R})$-strand dynamics admit singular solutions associated with a~pair of momentum maps.
Then we shall derive the one-peakon solution and solve for the collision dynamics of two-peakon interactions on
$G$-strands.
In Section~\ref{complex-sec} we introduce the complexif\/ied ${\rm Dif\/f}(\mathbb{R})$-strand equations and solve
their partial dif\/ferential equations corresponding to two-peakon collision dynamics.
Thus, Section~\ref{complex-sec} investigates the interactions between singular solutions of the $G$-strand equations
for the group of dif\/feomorphisms ${\rm Dif\/f}(\mathbb{R})$, both as evolutionary equations and as harmonic maps.
Finally, Section~\ref{conclusion-sec} summarizes the present results and provides some outlook for future research.

\section[The $G$-strand PDE for $G={\rm Dif\/f}(\mathbb{R})$]{The $\boldsymbol{G}$-strand PDE for $\boldsymbol{G={\rm Dif\/f}(\mathbb{R})}$}
\label{DiffStrand-sec}

This section studies the $G$-strand system that arises when we choose $G={\rm Dif\/f}(\mathbb{R})$ and the Lagrangian
involves the $H^1$ Sobolev norm.
This case is reminiscent of f\/luid dynamics and may be written naturally in terms of right-invariant tangent vectors
$u(t,s,x)$ and ${v}(t,s,x)$ def\/ined~by
\begin{gather*}
\partial_t{g}={u}\circ g
\qquad
\hbox{and}
\qquad
\partial_s{g}={v}\circ g,
\end{gather*}
where the symbol $\circ$ denotes composition of functions.
In this right-invariant case, the $G$-strand system of partial dif\/ferential equations (PDE) with reduced Lagrangian
$\ell({u},{v})$ takes the fol\-lo\-wing form, which generalizes to $\mathbb{R}^d$ in any number of spatial dimensions,
\begin{gather}
\frac{\partial}{\partial t}\frac{\delta\ell}{\delta{u}}+\frac{\partial}{\partial s}\frac{\delta\ell}{\delta{v}} =
- {\rm ad}^*_{u}\frac{\delta\ell}{\delta{u}}-{\rm ad}^*_{{v}}\frac{\delta\ell}{\delta{v}},
\qquad
\frac{\partial{v}}{\partial t}-\frac{\partial{u}}{\partial s} ={\rm ad}_{u}{v}.
\label{Gstrand-eqn1R}
\end{gather}
The distinction between the maps $({u},{v}):\mathbb{R}\times\mathbb{R}\to\mathfrak{g}\times\mathfrak{g}$ and their
pointwise values $({u}(t,s),{v}(t,s))$ $\in\mathfrak{g}\times\mathfrak{g}$ will always be clear in context, so that no
confusion will arise.
Likewise, for the variational derivatives ${\delta\ell}/{\delta{{u}}}$ and ${\delta\ell}/{\delta{v}}$.

Equations~\eqref{Gstrand-eqn1R} form a~subset of the equations studied in~\cite{GBRa2008,Ho2002,HoKu1988,Tr2012} for
complex f\/luids.
They also form a~subset of the equations for molecular strands studied in~\cite{ElGBHoPuRa2010}.
The latter comparison further justif\/ies the name \emph{$G$-strands} for the systems being studied here.
These equations may be derived from Hamilton's principle for an af\/f\/ine Lie group action, under which the auxiliary
equation for ${v}$ may be interpreted as an advection law.
This interpretation is discussed further in~\cite{ElGBHoPuRa2010,GBRa2008}.
The covariant form of these equations has been developed in~\cite{FGB2012}.

\subsection[The $G$-strand Hamiltonian structure]{The $\boldsymbol{G}$-strand Hamiltonian structure}

Upon setting $m={\delta\ell}/{\delta{u}}$ and
$n={\delta\ell}/{\delta{v}}$, the right-invariant $G$-strand equations in~\eqref{Gstrand-eqn1R} for maps
$\mathbb{R}\times\mathbb{R}\to G={\rm Dif\/f}(\mathbb{R})$ in one spatial dimension may be expressed as a~system of two
1+2 PDEs in $(t,s,x)$,
\begin{gather}
m_t+n_s =- {\rm ad}^*_{u}m-{\rm ad}^*_{v}n=-({u}m)_x-m{u}_x-({v}n)_x-n{v}_x,
\nonumber \\
{v}_t-{u}_s =- {\rm ad}_{v}{u}=-{u}{v}_x+{v}{u}_x.
\label{Gstrand-eqn2R}
\end{gather}
The corresponding Hamiltonian structure for these ${\rm Dif\/f}(\mathbb{R})$-strand equations is obtained by Legendre
transforming~\cite{HoMaRa1998} to
\begin{gather*}
h(m,{v})=\langle m, {u}\rangle-\ell({u}, {v}).
\end{gather*}
One may then write the $m$-${v}$ equations~\eqref{Gstrand-eqn2R} in Lie--Poisson Hamiltonian form as
\begin{gather}
\frac{d}{dt}
\begin{bmatrix}
m
\\
{v}
\end{bmatrix}
=
\begin{bmatrix}
- {\rm ad}^*_\Box m&\partial_s+{\rm ad}^*_{v}
\\
\partial_s-{\rm ad}_{v}&0
\end{bmatrix}
\begin{bmatrix}
{\delta h}/{\delta m}={u}
\\
{\delta h}/{\delta{v}}=- n
\end{bmatrix},
\label{1stHamForm}
\end{gather}
where $({\rm ad}^*_\Box m)v:={\rm ad}^*_v m$ for $v\in\mathfrak{g}$ and $m\in\mathfrak{g}^*$. This is the Lie--Poisson
bracket dual to the action of the semidirect-product Lie algebra
\begin{gather*}
\mathfrak{g}=\mathfrak{X}(\mathbb{R})\,\circledS\,\Lambda^1({\rm Dens})(\mathbb{R})\oplus C(\partial_s),
\end{gather*}
in which $\mathfrak{X}(\mathbb{R})$ is the space of vector f\/ields and $\Lambda^1({\rm Dens})(\mathbb{R})$ is the
space of 1-form densities on the real line $\mathbb{R}$, plus a~generalized 2-cocycle $C(\partial_s)$.
In the lower of\/f-diagonal entry of the Hamiltonian matrix in~\eqref{1stHamForm}, one recognizes the vector-f\/ield
covariant derivative in $s$, and f\/inds its negative adjoint operator in the upper of\/f-diagonal entry.
The Lie--Poisson bracket in~\eqref{1stHamForm} f\/irst arose in the Hamiltonian formulation of chromohydrodynamics,
i.e., the dynamics of a~Yang--Mills f\/luid plasma~\cite{GiHoKu1982,GiHoKu1983}.
For discussions of how such Lie--Poisson Hamiltonian structures arise in complex f\/luids with f\/inite-dimensional
order-parameter groups (broken symmetries) and full discussions of their Lie algebraic properties,
see~\cite{GBRa2008,Ho2002,HoKu1988,Tr2012} and references therein.

{\bf Relation to the Camassa--Holm equation~\cite{CaHo1993}.} An {interesting subcase} of the system of
semi-stationary ${\rm Dif\/f}(\mathbb{R})$-strand equations~\eqref{Gstrand-eqn1R} arises when one chooses the
Lagrangian $\ell({u},{v})$ in~\eqref{Gstrand-eqn1R} to depend only on ${u}$, as its $H^1$ norm on the real line,
\begin{gather*}
\ell({u},{v})=\frac12\|{u}\|^2_{H^1},
%\label{CH-pkn-Lag}
\end{gather*}
with vanishing boundary conditions, as $|x|\to\infty$.
In this case, $m={u}-{u}_{xx}$, and this restriction of the equations in~\eqref{Gstrand-eqn1R} provides an extension of
the completely integrable Camassa--Holm (CH) equation~\cite{CaHo1993},
\begin{gather}
m_t=- {\rm ad}^*_{{u}}m=-({u}m)_x-m{u}_x
\qquad
\hbox{with}
\qquad
m=\frac{\delta\ell}{\delta{u}}={u}-{u}_{xx}.
\label{CH-eqn}
\end{gather}
Specif\/ically, these modif\/ied $G$-strand equations reduce in the absence of ${v}$-dependence precisely to the CH
equation, which admits singular solutions known as \emph{peakons} in the form~\cite{CaHo1993}
\begin{gather}
m(x,t)=\sum_a M_a(t)\delta\big(x-Q^a(t)\big),
\label{pkn-soln}
\end{gather}
where we sum in $a\in\mathbb{Z}$ over the integers, or over any subset of the integers.
The peakon solution~\eqref{pkn-soln} of the CH equation may be understood as a~singular momentum map obtained from the
left action of ${\rm Dif\/f}(\mathbb{R})$ on embeddings of points on the real line~\cite{HoMa2004}.
The generalization to higher dimensions is also possible.
See~\cite{HoMa2004} for further discussion.

\subsection[Peakon solutions of the ${\rm Dif\/f}(\mathbb{R})$-strand equations in~(\ref{Gstrand-eqn2R})]{Peakon solutions of the $\boldsymbol{{\rm Dif\/f}(\mathbb{R})}$-strand equations in~(\ref{Gstrand-eqn2R})}

 With the
following choice of Lagrangian using the $H^1$ norm,
\begin{gather}
\ell(u,v)=\frac12\|u\|^2_{H^1}+\frac{\sigma}{2}\|v\|^2_{H^1},
\label{Gstrand-pkn-Lag}
\end{gather}
the ${\rm Dif\/f}(\mathbb{R})$-strand equations~\eqref{Gstrand-eqn2R} admit peakon solutions in \emph{both} momenta $m$
and $n$, with continuous velocities $u$ and $v$.
Here $\sigma=\pm1$ provides two possibilities for a~sign choice.
Note that only the choice $\sigma=-1$ leads to a~positive-def\/inite Hamiltonian.
The choice $\sigma=1$ in~\eqref{Gstrand-pkn-Lag} yields a~conserved Hamiltonian that is not bounded below.
We state the main result in the following theorem.
\begin{theorem}\label{Gstrand-HP}
The ${\rm Dif\/f}(\mathbb{R})$-strand equations~\eqref{Gstrand-eqn2R} admit singular solutions expressible as linear
superpositions summed over $a\in\mathbb{Z}$~{\rm \cite{HoIvPe2012}}
\begin{gather}
m(x,t,s)=\sum_a M_a(t,s)\delta\big(x-Q^a(t,s)\big),
\nonumber\\
n(x,t,s) =\sum_a N_a(t,s)\delta\big(x-Q^a(t,s)\big),
\nonumber\\
u(x,t,s) =K*m=\sum_a M_a(t,s)K\big(x,Q^a\big),
\nonumber\\
v(x,t,s) =\sigma K*n=\sigma\sum_a N_a(t,s)K\big(x,Q^a\big),
\label{Gstrand-singsolns}
\end{gather}
that are \emph{peakons} in the case that $K(x,y)=\frac12e^{-|x-y|}$ is the Green function of the Helmholtz operator
$1-\partial_x^2$.
These singular solutions follow from Hamilton's principle $\delta S=0$ for the constrained action $S=\int L(u,v,Q)\,dt$
given by
\begin{gather*}
S=
\int\ell(u,v)+\sum_a M_a(t,s)\big(\partial_t Q^a(t,s)-u(Q^a,t,s)\big)\\
\hphantom{S=}{}
+\sum_a N_a(t,s)\big(\partial_s Q^a(t,s)-v(Q^a,t,s)\big)\,dsdt.
\end{gather*}
and they are \emph{peakons} in $(t,s)$ for the Lagrangian $\ell(u,v)$ given in equation~\eqref{Gstrand-pkn-Lag}.
\end{theorem}

\begin{proof}
The proof of this theorem is a~direct calculation following the method of momentum maps and the Clebsch approach for
this type of problem introduced in~\cite{HoIvPe2012}.
The results for the present case are given below.
\end{proof}

The solution parameters $\{Q^a(t,s),M_a(t,s),N_a(t,s)\}$ with $a\in\mathbb{Z}$ that specify the singular
solutions~\eqref{Gstrand-singsolns} are determined by the following set of evolutionary PDEs in~$s$ and~$t$, in which
we denote $K^{ab}:=K(Q^a,Q^b)$ with integer summation indices $b,c,e\in\mathbb{Z}$:
\begin{gather}
\partial_t Q^a(t,s) =u(Q^a,t,s)=\sum_b M_b(t,s)K^{ab},
\nonumber\\
\partial_s Q^a(t,s) =v(Q^a,t,s)=\sigma\sum_b N_b(t,s)K^{ab},
\nonumber \\
\partial_t M_a(t,s) =- \partial_s N_a-\sum_c(M_aM_c+\sigma N_aN_c)\frac{\partial K^{ac}}{\partial Q^a}
\qquad
\hbox{(no sum on~$a$)},
\nonumber\\
\partial_t N_a(t,s) =
\sigma\partial_s M_a+\sum_{b,c,e}(N_bM_c-M_bN_c)\frac{\partial K^{ec}}{\partial Q^e}\big(K^{eb}-K^{cb}\big)\big(K^{-1}\big)_{ae}.
\label{Gstrand-eqns}
\end{gather}

The notation $(K^{-1})_{ab}$ is for the $ab$-entry of the matrix $K^{-1}$.

The last pair of equations in~\eqref{Gstrand-eqns} may be solved as a~system for the momenta, or Lagrange multipliers
$(M_a,N_a)$, then used in the previous pair to update the support set of positions $Q^a(t,s)$.
Given $Q^a(t,s)$ for $a\in\mathbb{Z}$, one constructs $(m(x,t,s),n(x,t,s))$ along the solution paths $x=Q^a(t,s)$ from
the f\/irst pair of~\eqref{Gstrand-singsolns} and then obtains $(u(x,t,s),v(x,t,s))$ for $x\in\mathbb{R}$ from the
second pair.
Alternatively, knowing the position $Q^a(t,s)$, $a\in\mathbb{Z}$, for all $s$ at a~given time $t$, also determines
$N_a$ upon inverting the matrix $K^{ab}$ in the second equation in~\eqref{Gstrand-eqns}.
The values of~$M_a$ can be determined similarly.

{\bf Single-peakon solution.} The single-peakon solution of the ${\rm Dif\/f}(\mathbb{R})$-strand
equations~\eqref{Gstrand-eqn2R} is straightforward to obtain from~\eqref{Gstrand-eqns}.
Namely, writing the f\/irst two equations for a~single peakon and substituting into the third equation shows that
$Q^1(t,s)$ satisf\/ies the equation,
\begin{gather*}
\big(\partial_s^2+\sigma\partial_t^2\big)Q^1(t,s)=0.
%\label{Q-1-peak}
\end{gather*}
Thus, any solution $h(t,s)$ of the linear equation $(\partial_s^2+\sigma\partial_t^2)h(t,s)=0$ provides a~solution
$Q^1=h(t,s)$.
One then f\/inds from the f\/irst two equations in~\eqref{Gstrand-eqns} that
\begin{gather*}
M^1(t,s)=\frac{1}{K_0}h_t(t,s),
\qquad
N^1(t,s)=\frac{\sigma}{K_0}h_s(t,s),
%\label{MN-1-peak}
\end{gather*}
where $K_0=K(0)$ with $K(X-Y)\equiv K(X,Y)$.
The solution for the single-peakon parameters $Q^1$, $M^1$ and $N^1$ depends only on the choice of the harmonic function
$h(t,s)$, which in turn depends on the~$(t,s)$ boundary conditions.
In particular, the single-peakon does not depend on the shape of the Green's function~$K(x)$.

\subsection[Two-peakon interactions on $G$-strands]{Two-peakon interactions on $\boldsymbol{G}$-strands}
\label{2peakon-collision}

Consider an operator norm $\|w\|^2_{\mathcal{Q}}=\frac12\int w \mathcal{Q} w\,dx$ where the operator $\mathcal{Q}$ is
symmetric, positive-def\/inite and spatially translation-invariant, with vanishing boundary conditions, as
$|x|\to\infty$.
For example, the $H^1$ norm for the CH peakons in the Lagrangian~\eqref{Gstrand-pkn-Lag} arises when $\mathcal{Q}$ is
chosen to be the Helmholtz operator, so that $\mathcal{Q}w=w-w_{xx}$.

Let the Lagrangian be
\begin{gather*}
\ell(u,v)=\frac12\|u\|^2_{\mathcal{Q}}+\frac{\sigma}{2}\|v\|^2_{\mathcal{Q}}.
%\label{Q-op-Lag}
\end{gather*}

Consider peakons at positions $Q^1(t,s)$ and $Q^2(t,s)$ on the real line.
The Green's function $K$ for the operator $\mathcal{Q}$ depends only on the dif\/ference of the peakon positions,
$X(t,s)=Q^1-Q^2$.
Then the f\/irst two equations in~\eqref{Gstrand-eqns} imply
\begin{gather}
\partial_t X =(M_1-M_2)(K_0-K(X)),
\qquad
\partial_s X =\sigma(N_1-N_2)(K_0-K(X)),
\label{Qdiff-eqns}
\end{gather}
where $K_0=K(0)$ is the value of the Green's function $K(X)$ for the operator $\mathcal{Q}$ when the peakon positions
coincide, so that $X=0$.

The second pair of equations in~\eqref{Gstrand-eqns} may then be written as
\begin{gather}
  \partial_t M_1 =-\partial_s N_1-(M_1M_2+\sigma N_1N_2)K'(X),
\nonumber\\
\partial_t M_2 =-\partial_s N_2+(M_1M_2+\sigma N_1N_2)K'(X),
\nonumber\\
  \partial_t N_1 =\sigma\partial_s M_1+(N_1M_2-M_1N_2)\frac{K_0-K}{K_0+K}K'(X),
\nonumber\\
\partial_t N_2 =\sigma\partial_s M_2+(N_1M_2-M_1N_2)\frac{K_0-K}{K_0+K}K'(X).
\label{Gstrand-pp}
\end{gather}
Asymptotically, when the peakons are far apart, the system~\eqref{Gstrand-pp} simplif\/ies, since
$\frac{K_0-K}{K_0+K}\to1$ and $K'(X)\to0$ as $|X|\to\infty$.

The system~\eqref{Gstrand-pp} has two immediate conservation laws obtained from their sums and dif\/fe\-ren\-ces,
\begin{gather*}
\partial_t(M_1+M_2) =- \partial_s(N_1+N_2),
\qquad
\partial_t(N_1-N_2) =\sigma\partial_s(M_1-M_2).
%\label{pp-CLs}
\end{gather*}
These may be resolved by setting
\begin{alignat}{3}
& M_1-M_2 =\frac{\partial_t X}{K_0-K},
\qquad &&
N_1-N_2=\sigma\frac{\partial_s X}{K_0-K}, &
\nonumber\\
& M_1+M_2=\partial_s\phi,
\qquad &&
N_1+N_2=- \partial_t\phi, &
\label{pp-Xpotentials}
\end{alignat}
and introducing two potential functions, $X$ and $\phi$, for which equality of cross derivatives will now produce the
system of equations~\eqref{Qdiff-eqns} and~\eqref{Gstrand-pp}.
Namely, from~\eqref{pp-Xpotentials}, one may solve for
\begin{alignat}{3}
& M_1 =\frac{\partial_t X}{2(K_0-K)}+\frac12\partial_s\phi,
\qquad &&
N_1=\frac{\sigma\partial_s X}{2(K_0-K)}-\frac12\partial_t\phi,&
\nonumber\\
& M_2 =- \frac{\partial_t X}{2(K_0-K)}+\frac12\partial_s\phi,
\qquad &&
N_2=- \frac{\sigma\partial_s X}{2(K_0-K)}-\frac12\partial_t\phi. &
\label{pp-Xpotentials1}
\end{alignat}
Now substituting into the f\/irst and third equations of the system~\eqref{Gstrand-pp} yields the determining equations
for the potentials $X$ and $\phi$,{\samepage
\begin{gather}
\big(\partial_t^2+\sigma\partial_s^2\big)\phi +\frac{K'}{K_0+K}(X_t\phi_t+\sigma X_s\phi_s)=0,
\nonumber\\
\big(\partial_t^2+\sigma\partial_s^2\big)X +\frac{K'}{2(K_0-K)}\big(X_t^2+\sigma X_s^2\big)=
- \frac12K'(K_0-K)\big(\phi_s^2+\sigma\phi_t^2\big),
\label{pp-Xpotentials2}
\end{gather}}

{\bf A simplif\/ication.} A simplif\/ication arises if $\phi=0$, in which case the collision is perfectly
antisymmetric, as seen from equation~\eqref{pp-Xpotentials}.
This is the peakon-antipeakon collision, for which the second equation in~\eqref{pp-Xpotentials2} reduces to
\begin{gather*}
\big(\partial_t^2+\sigma\partial_s^2\big)X+\frac{K'}{2(K_0-K)}\big(X_t^2+\sigma X_s^2\big)=0.
\end{gather*}
This equation can be easily rearranged to produce a~linear equation:
\begin{gather*}
\big(\partial_t^2+\sigma\partial_s^2\big)F(X)=0,
\qquad
\hbox{where}
\qquad
F(X)=\int_{X_0}^X(K_0-K(Y))^{-1/2}\,dY.
\end{gather*}
Thus, the dynamics of the relative spacing $X(t,s)=Q^1-Q^2$ in the peakon-antipeakon collision may be obtained by
elementary means.
In particular, for the $H^1$ peakon case when $K(Y)=\frac{1}{2}e^{-|Y|}$, we have
\begin{gather}
F(X)=\sqrt{2}\int_{X_0}^X\frac{1}{\sqrt{1-e^{-|Y|}}}\,dY.
\label{F-express}
\end{gather}

We can take for simplicity $X_0=0$, this would change $F(X)$ only by a~constant.
When $X>0$ the integral becomes
\begin{gather*}
F(X)=\sqrt{2}\int_{0}^X\frac{1}{\sqrt{1-e^{-Y}}}\,dY=\sqrt{2}\int_{0}^X\frac{e^{Y/2}}{\sqrt{e^{Y}-1}}\,dY
\\
\hphantom{F(X)}{} =2\sqrt{2}\int_{0}^X\frac{d e^{Y/2}}{\sqrt{(e^{Y/2})^2-1}}=2\sqrt{2}\cosh^{-1}\big(e^{X/2}\big).
\end{gather*}

In general $F(X)=2\sqrt{2}\,\text{sign}(X)\cosh^{-1}\left(e^{|X|/2}\right)$.
Hence the solution $X(t,s)$ can be expressed in terms of any given solution $h(t,s)$ of the linear equation
$(\partial_t^2+\sigma\partial_s^2)h(t,s)=0$ as
\begin{gather*}
X(t,s)=\pm\ln\left({\rm cosh}^2(h(t,s))\right).
%\label{X-express}
\end{gather*}
When $\sigma=1$, $h(t,s)$ is any harmonic function, when $\sigma=-1$, $h(t,s)$ is any solution of the wave equation.

\section[Complexif\/ied $G$-strand ${\rm Dif\/f}(\mathbb{R})$ equations]{Complexif\/ied $\boldsymbol{G}$-strand $\boldsymbol{{\rm Dif\/f}(\mathbb{R})}$ equations}
\label{complex-sec}

As mentioned in the Introduction, the class of problems being treated here generalizes the classical chiral model,
which is a~special case of the $G$-strands when their structure group is ${\rm SO}(3)$, as discussed in~\cite{HoIvPe2012}.
Harmonic maps are also related to the chiral model as discussed in~\cite{Uh1989}, upon complexifying the independent
coordinates of the chiral model.
Pursuing this idea, we shall complexify the independent coordinates $(t,s)$ of the $G$-strand
equations~\eqref{Gstrand-eqn2R} on the Lie group of dif\/feomorphisms ${\rm Dif\/f}(\mathbb{R})$ and thereby derive
a~class of equations which may be regarded as partial dif\/ferential equations for harmonic maps of $\mathbb{C}$ into
${\rm Dif\/f}(\mathbb{R})$.
This is potentially a~rich area for further mathematical study.

\subsection[Complexifying coordinates $(t,s)$ for a~real Lagrangian]{Complexifying coordinates $\boldsymbol{(t,s)}$ for a~real Lagrangian}

We complexify coordinates
$(t,s)\in\mathbb{R}^2\to(z,\bar{z})\in\mathbb{C}$, where $\bar{z}$ denotes the complex conjugate of~$z$ and identify
$v=\bar{u}$.
Consequently, the Euler--Poincar\'e $G$-strand equations in~\eqref{Gstrand-eqn2R} become
\begin{gather}
\frac{\partial}{\partial z}\frac{\delta\ell}{\delta u}+\frac{\partial}{\partial\bar{z}}\frac{\delta\ell}{\delta\bar{u}}=
-\,{\rm ad}^*_{u}\frac{\delta\ell}{\delta u}-{\rm ad}^*_{\bar{u}}\frac{\delta\ell}{\delta\bar{u}},
\qquad
\hbox{and}
\qquad
\frac{\partial\bar{u}}{\partial z}-\frac{\partial u}{\partial\bar{z}}={\rm ad}_u\bar{u}.
\label{ComplexGstrand-eqns1}
\end{gather}

Here we assume that the Lagrangian $\ell$ is def\/ined as
\begin{gather}
\ell(u,\bar{u})=\frac12\|u\|_\mathcal{Q}^2+\frac{\sigma}{2}\|\bar{u}\|_\mathcal{Q}^2,
\label{ComplexNorm-ell}
\end{gather}
where the operator $\mathcal{Q}$ and the expression $\|w\|_\mathcal{Q}$ are formally def\/ined as before, although
$\|w\|_\mathcal{Q}$ does not def\/ine a~norm for complex $w$.

By evaluating $m={\delta\ell}/{\delta u}=\mathcal{Q}u$,
$\bar{m}={\delta\ell}/{\delta\bar{u}}=\sigma\mathcal{Q}\bar{u}$, we observe that consistency of this construction
requires that we set $\sigma=1$ in equation~\eqref{ComplexNorm-ell}.

For the real Lagrangian $\ell$, equations~\eqref{ComplexGstrand-eqns1} may be rewritten as
\begin{gather}
m_{z}+\bar{m}_{\bar{z}} =- {\rm ad}^*_{u}m-{\rm ad}^*_{\bar{u}}\bar{m}=
-(u m)_x-mu_x-(\bar{u}\bar{m})_x-\bar{m}\bar{u}_x,
\nonumber\\
\bar{u}_z-u_{\bar{z}} =- {\rm ad}_{\bar{u}}u=-u\bar{u}_x+\bar{u}u_x,
\label{ComplexGstrand-eqns2}
\end{gather}
where the independent coordinate $x\in\mathbb{R}$ is on the real line, although coordinates $(z,\bar{z})\in\mathbb{C}$
are complex, as are the solutions $u$, and $m=\mathcal{Q}u$, of the equations~\eqref{ComplexGstrand-eqns1} for the real
Lagrangian~\eqref{ComplexNorm-ell} for $\sigma=1$.
\begin{remark}\qquad
\begin{itemize}\itemsep=0pt
\item Equations~\eqref{ComplexGstrand-eqns2} are invariant under two involutions, $P$ and $C$, where
\begin{gather*}
P: \ (x,m)\to(-x,-m)
\qquad
\hbox{and}
\qquad
C: \ \hbox{complex conjugation.}
\end{gather*}
\item For real variables $m=\bar{m}$, $u=\bar{u}$ and real evolution parameter $z=\bar{z}=:t$,
equations~\eqref{ComplexGstrand-eqns2} reduce to the Camassa--Holm equation~\eqref{CH-eqn}.
\end{itemize}
\end{remark}

\subsection[Peakon solutions of complexif\/ied ${\rm Dif\/f}(\mathbb{R})$-strand equations
in~(\ref{ComplexGstrand-eqns2})]{Peakon solutions of complexif\/ied $\boldsymbol{{\rm Dif\/f}(\mathbb{R})}$-strand equations
in~(\ref{ComplexGstrand-eqns2})}

The complexif\/ied ${\rm Dif\/f}(\mathbb{R})$-strand equations~\eqref{ComplexGstrand-eqns2} admit peakon solutions in
complex momenta~$m$, with continuous complex velocity $u$.
We state this result in the following theorem.
\begin{theorem}\label{Gstrand-HP+}
The complexified ${\rm Dif\/f}(\mathbb{R})$-strand equations~\eqref{ComplexGstrand-eqns2} admit singular solutions
expressed as linear superpositions summed over $a\in\mathbb{Z}$
\begin{gather}
m(x,z,\bar{z})=\sum_a M_a(z,\bar{z})\delta\big(x-Q^a(z,\bar{z})\big),
\nonumber \\
u(x,z,\bar{z}) =K*m=\sum_a M_a(z,\bar{z})K\big(x,Q^a\big).
\label{Gstrand-singsolns2}
\end{gather}
These are peakon solutions for $K(x,y)=\frac12e^{-|x-y|}$.
An important feature of these singular solutions is that $Q^a$ with $a\in\mathbb{Z}$ are \emph{real} functions of $z$
and $\bar{z}$.
These singular solutions follow from Hamilton's principle $\delta S=0$ for the constrained action $S=\int
L(u,\bar{u},Q)\,dt$ given by
\begin{gather*}
S=
\int\ell(u,\bar{u})
+\Re\left(\sum_a M_a(z,\bar{z})\big(\partial_z Q^a(z,\bar{z})-u(Q^a,z,\bar{z})\big)\right)dz d\bar{z},
\end{gather*}
where $\Re$ is the real part of the expression.
\end{theorem}

\begin{proof}
The proof follows the same lines as the proof in~\cite{HoIvPe2012} for the real case, keeping in
mind that the coordinate $x\in\mathbb{R}$ is a~real independent variable.

The solution parameters $Q^a(z,\bar{z})\in\mathbb{R}$ and $M_a(z,\bar{z})\in\mathbb{C}$ with $a\in\mathbb{Z}$ that
specify the singular solutions~\eqref{Gstrand-singsolns2} are determined by the following set of evolutionary PDEs in
$z$ and $\bar{z}$, in which we denote $K^{ab}:=K(Q^a,Q^b)$ with integer summation indices $b,c,e\in\mathbb{Z}$:
\begin{gather}
\partial_z Q^a(z,\bar{z}) =u(Q^a,z,\bar{z})=\sum_b M_b(z,\bar{z})K^{ab},
\nonumber\\
\partial_z M_a(z,\bar{z}) =
- \partial_{\bar{z}}\bar{M}_a-\sum_c(M_aM_c+\bar{M}_a\bar{M}_c)\frac{\partial K^{ac}}{\partial Q^a}
\qquad
\hbox{(no sum on~$a$),}
\nonumber\\
\partial_z\bar{M}_a(z,\bar{z}) =
\partial_{\bar{z}}M_a+\sum_{b,c,e}(\bar{M}_bM_c-M_b\bar{M}_c)\frac{\partial K^{ec}}{\partial Q^e}\big(K^{eb}-K^{cb}\big)\big(K^{-1}\big)_{ae}.
\label{Gstrand-eqns2}
\end{gather}

The last pair of equations in~\eqref{Gstrand-eqns2} may be solved as a~system for the momenta, or Lagrange multipliers
$(M_a)$, then used in the f\/irst equation to update the support set of positions $Q^a(z,\bar{z})$.
Given $Q^a(z,\bar{z})$ for $a\in\mathbb{Z}$, one constructs $m(x,z,\bar{z})$ along the solution paths
$x=Q^a(z,\bar{z})$ from the f\/irst equation of~\eqref{Gstrand-singsolns2} and then obtains $u(x,z,\bar{z})$ for
$x\in\mathbb{R}$ from the second equation of~\eqref{Gstrand-singsolns2}.
Alternatively, knowing the position $Q^a(z,\bar{z})$, $a\in\mathbb{Z}$, for all $z$ and $\bar{z}$, one also determines
$\bar{M}_a$ upon inverting the matrix $K^{ab}$ in the third equation in~\eqref{Gstrand-eqns2}.
Note that the third equation in~\eqref{Gstrand-eqns2} is not the complex conjugate of the previous one.
\end{proof}

{\bf Single-peakon solution.} The single-peakon solution is again straightforward to obtain from
equations~\eqref{Gstrand-eqns2}.
The solution $Q^1(z,\bar{z})$ satisf\/ies the (linear) two-dimensional wave equation.

{\bf Two-peakon solution.} The two-peakon collision for the complexif\/ied $G$-strand equation emerges by
following the same path as in the previous section~\ref{2peakon-collision} with appropriate minor changes.
In particular, equation~\eqref{pp-Xpotentials1} becomes
\begin{gather*}
M_1 =\frac{\partial_z X}{2(K_0-K)}+\frac{i}{2}\partial_{\bar{z}}\phi,
\qquad
M_2 =- \frac{\partial_z X}{2(K_0-K)}+\frac{i}{2}\partial_{\bar{z}}\phi,
%\label{pp-Xpotentials1a}
\end{gather*}
where $\phi(z,\bar{z})$ is real.
Now substituting into the second and third equations of the system~\eqref{Gstrand-eqns2} we obtain
\begin{gather*}
\phi_{zz}+\phi_{\bar{z}\bar{z}} +\frac{K'}{K_0+K}(X_z\phi_z+X_{\bar{z}}\phi_{\bar{z}})=0,
\\
X_{zz}+X_{\bar{z}\bar{z}} +\frac{K'}{2(K_0-K)}\big(X_z^2+X_{\bar{z}}^2\big)=\frac12K'(K_0-K)\big(\phi_z^2+\phi_{\bar{z}}^2\big).
%\label{pp-Xpotentials2a0}
\end{gather*}
Introducing the notation
\begin{gather*}
\partial_z=\partial_\xi+i\partial_\eta
\qquad
\hbox{and}
\qquad
\partial_{\bar{z}}=\partial_\xi-i\partial_\eta
\end{gather*}
we obtain the following PDE system for the real variables $\xi$, $\eta$:
\begin{gather*}
\phi_{\xi\xi}-\phi_{\eta\eta} +\frac{K'}{K_0+K}(X_\xi\phi_\xi-X_\eta\phi_\eta)=0,
\\
X_{\xi\xi}-X_{\eta\eta} +\frac{K'}{2(K_0-K)}\big(X_\xi^2-X_\eta^2\big)=\frac12K'(K_0-K)\big(\phi_\xi^2-\phi_\eta^2\big).
%\label{pp-Xpotentials2a}
\end{gather*}
When $\phi=0$, we again have peakon-antipeakon collisions with a~simple solution as in Section~\ref{2peakon-collision}.
The equation reduces to the two-dimensional (linear) wave equation:
\begin{gather}
\left(\partial_\xi^2-\partial_\eta^2\right)F(X)=0,
\label{wave-soln}
\end{gather}
where $F$ is the same function as in~\eqref{F-express}.
The obvious dif\/ference when $\sigma=1$ is that the Laplacian $\partial_s^2+\partial_t^2$ in~\eqref{pp-Xpotentials2}
for real parameters $(t,s)$ is replaced in~\eqref{wave-soln} by the wave operator $\partial_\xi^2-\partial_\eta^2$ for
this complex case.

\subsection[Complexif\/ied ${\rm Dif\/f}(\mathbb{R})$-strand equations for a~real Lagrangian given by a~complex norm]{Complexif\/ied $\boldsymbol{{\rm Dif\/f}(\mathbb{R})}$-strand equations for a~real Lagrangian\\ given by a~complex norm}

Here the Lagrangian $\ell$ is taken to be real, and def\/ined by a~complex norm as
\begin{gather}\label{ComplexNorm-ell 1}
\ell(u,\bar{u})=\int u\mathcal{Q}\bar{u}\,dx,
\end{gather}
where the operator $\mathcal{Q}$ is def\/ined as before -- a~real, symmetric, positive-def\/inite operator, invariant
under spatial translations in $x\in\mathbb{R}$, and with vanishing boundary conditions as $|x|\to\infty$.

Since the Lagrangian $\ell$ is chosen to be a~norm on complex functions def\/ined by the spatial dif\/ferential
operator $\mathcal{Q}$, we may regard the f\/irst equation in the complexif\/ied $G$-strand
system~\eqref{ComplexGstrand-eqns1} as describing geodesic motion according to the map $(z,\bar{z})\in\mathbb{C}\to{\rm
Dif\/f}(\mathbb{R})$ with respect to the above norm.
The second equation~\eqref{ComplexGstrand-eqns1} may be regarded as a~compatibility (or zero-curvature) condition in
the complex parameters $z$ and $\bar{z}$.

Evaluating ${\delta\ell}/{\delta u}=\mathcal{Q}\bar{u}$, ${\delta\ell}/{\delta\bar{u}}=\mathcal{Q}u$, the
equations~\eqref{ComplexGstrand-eqns1} for the real Lagrangian~\eqref{ComplexNorm-ell 1} may be rewritten as
\begin{gather}
\bar{\mu}_{z}+\mu_{\bar{z}} =- {\rm ad}^*_{u}\bar{\mu}-{\rm ad}^*_{\bar{u}}\mu=
-(u\bar{\mu})_x-\bar{\mu}u_x-(\bar{u}\mu)_x-\mu\bar{u}_x,
\nonumber\\
\bar{u}_z-u_{\bar{z}} =- {\rm ad}_{\bar{u}}u=-u\bar{u}_x+\bar{u}u_x,
\label{ComplexGstrand-eqns21}
\end{gather}
where $\mu=\mathcal{Q}u$ for convenience.
\begin{remark}\qquad
\begin{itemize}\itemsep=0pt
\item Equations~\eqref{ComplexGstrand-eqns21} are invariant under the involutions, $P$ and $C$, where
\begin{gather*}
P: \ (x,\mu)\to(-x,-\mu)
\qquad
\hbox{and}
\qquad
C: \ \hbox{complex conjugation.}
\end{gather*}
\item For real variables $\mu=\bar{\mu}:=m$, $u=\bar{u}$ and real evolution parameter $z=\bar{z}=:t$,
equations~\eqref{ComplexGstrand-eqns21} reduce to the Camassa--Holm equation~\eqref{CH-eqn}.
\end{itemize}
\end{remark}

Similarly to the previous case, the singular solutions expressed as linear superpositions summed over $a\in\mathbb{Z}$
\begin{gather*}
\mu(x,z,\bar{z}) =\sum_a M_a(z,\bar{z})\delta\big(x-Q^a(z,\bar{z})\big),
\\
u(x,z,\bar{z}) =K*\mu=\sum_a M_a(z,\bar{z})K\big(x,Q^a\big).
%\label{Gstrand-singsolns21}
\end{gather*}

The corresponding equations for the solution parameters $Q^a(z,\bar{z})\in\mathbb{R}$ and $M_a(z,\bar{z})\in\mathbb{C}$
with $a\in\mathbb{Z}$ are
\begin{gather}
\partial_z Q^a(z,\bar{z}) =u(Q^a,z,\bar{z})=\sum_b M_b(z,\bar{z})K^{ab},
\nonumber\\
\partial_z\bar{M}_a(z,\bar{z}) =
- \partial_{\bar{z}}M_a-\sum_c(\bar{M}_aM_c+M_a\bar{M}_c)\frac{\partial K^{ac}}{\partial Q^a}
\qquad
\hbox{(no sum on~$a$),}
\nonumber\\
\partial_z\bar{M}_a(z,\bar{z}) =
\partial_{\bar{z}}M_a+\sum_{b,c,e}(\bar{M}_bM_c-M_b\bar{M}_c)\frac{\partial K^{ec}}{\partial Q^e}\big(K^{eb}-K^{cb}\big)\big(K^{-1}\big)_{ae}.
\label{Gstrand-eqns21}
\end{gather}

The dif\/ference with the previous case is only in the second equation, since the compatibility equations are the same.

{\bf One-peakon solution.} For the one-peakon solution $\partial_{\bar{z}}M_1=0$, and therefore
$\partial_z\partial_{\bar{z}}Q^1=0$.
Thus $Q^1$ satisf\/ies the Laplace equation $\Delta Q^1=0$, and the solution is given by any harmonic function.

{\bf Two-peakon solution.} The two-peakon collision for this complexif\/ied $G$-strand equation can be obtained
in a~similar fashion.
This time the expressions are
\begin{gather*}
M_1 =\frac{\partial_z X}{2(K_0-K)}+\frac{i}{2}\partial_{z}\phi,
\qquad
M_2 =- \frac{\partial_z X}{2(K_0-K)}+\frac{i}{2}\partial_{z}\phi,
%\label{pp-Xpotentials1a1}
\end{gather*}
where $\phi(z,\bar{z})$ is real.
Now substituting into the second and third equations of the system~\eqref{Gstrand-eqns21} we obtain
\begin{gather*}
\begin{split}
\phi_{z\bar{z}}&+\frac{K'}{2(K_0+K)}(X_z\phi_{\bar{z}}+X_{\bar{z}}\phi_{z})=0,
\\ % may be a label omitted
X_{z\bar{z}}&+\frac{K'}{2(K_0-K)}X_z X_{\bar{z}}=-\frac12K'(K_0-K)\phi_z\phi_{\bar{z}}.
\end{split}
%\label{pp-Xpotentials2a01}
\end{gather*}
With the notations $\partial_z=\partial_\xi+i\partial_\eta$ and $\partial_{\bar{z}}=\partial_\xi-i\partial_\eta$ we
obtain the following system for the real variables $\xi$, $\eta$:
\begin{gather*}
\phi_{\xi\xi}+\phi_{\eta\eta} +\frac{K'}{K_0+K}(X_\xi\phi_\xi+X_\eta\phi_\eta)=0,
\\
X_{\xi\xi}+X_{\eta\eta} +\frac{K'}{2(K_0-K)}\big(X_\xi^2+X_\eta^2\big)=-\frac12K'(K_0-K)\big(\phi_\xi^2+\phi_\eta^2\big).
%\label{pp-Xpotentials2a1}
\end{gather*}
When $\phi=0$, we again have peakon-antipeakon collisions with a~simple solution as in Section~\ref{2peakon-collision}.
The equation reduces to the Laplace equation:
\begin{gather*}
\left(\partial_\xi^2+\partial_\eta^2\right)F(X)=0,
%\label{wave-soln}
\end{gather*}
where again $F$ is the same as in~\eqref{F-express}.

The role of the Laplacian has been studied before in connection with the relation between chiral models and harmonic
maps in~\cite{Uh1989}.
See~\cite{Gu2006} for further discussion and a~guide to the modern literature.
The corresponding studies for $G$-strands can be expected to be developed further in future work.
For example, it may be interesting to study the extension of $G$-strand equations to supersymmetric variables following
methods considered in~\cite{GuOl2006,Po2006}, particularly for $G={\rm Dif\/f}(\mathbb{R})$, if the singular solutions
would be preserved.

\section{Conclusions}
\label{conclusion-sec}

This paper has used the Euler--Poincar\'e (EP) framework for studying the $G$-strand equations, comprising the system of
partial dif\/ferential equations obtained from the EP variational equations~\eqref{Gstrand-eqn2R} for a~$G$-invariant
Lagrangian, coupled to an auxiliary \emph{zero-curvature} equation.
The latter has often been the departure point and main focus in other approaches, especially for the integrable chiral
models, where it sets up the Lax-pair formulation of the system~\cite{Lax68}.

Once the $G$-invariant Lagrangian has been specif\/ied, the combined dynamics and zero-curvature system of $G$-strand
equations in~\eqref{Gstrand-eqn2R} follows automatically in the EP framework.

The single-peakon solution of the ${\rm Dif\/f}(\mathbb{R})$-strand equations~\eqref{Gstrand-eqn2R} has been reduced to
a~solution of the linear Laplace or wave equations.
The stability of the single-peakon solution under perturbations into the full solution space of
equations~\eqref{Gstrand-eqn2R} would be an interesting problem for future work.

The pairwise dynamics of the singular peakon solutions for the $G$-strand equations with $G={\rm Dif\/f}(\mathbb{R})$
has also been studied.
The peakon-peakon and peakon-antipeakon collision interactions have been found to admit elementary solution methods
which reduce again to solving the linear Laplace or wave equations.

Finally, two complexif\/ied version of the ${\rm Dif\/f}(\mathbb{R})$-strand equations has been introduced and their
peakon collision solutions have been formulated and again solved by elementary means.
In the latter case, the analogous peakon-peakon and peakon-antipeakon collision interactions involve the solution of
the Laplace or wave operator in each case.

\subsection*{Acknowledgements}

We thank our friends A.M.~Bloch, C.J.~Cotter, F.~Gay-Balmaz, A.~Iserles, J.R.~Percival,
T.S.~Ratiu and C.~Tronci for their kind encouragement and thoughtful remarks during the course of this work.
We are thankful also to Dr.~Sergey Kushnarev and an anonymous referee whose comments and suggestions have helped us
a~lot in the revision of this paper.
DDH gratefully acknowledges partial support by the Royal Society of London's Wolfson Award scheme and the European
Research Council's Advanced Grant 267382 FCCA.
RII is supported by the Science Foundation Ireland (SFI), under Grant No.
09/RFP/MTH2144.

\pdfbookmark[1]{References}{ref}
\LastPageEnding

\end{document}